\documentclass[12pt,a4paper]{elsarticle}
\usepackage[T2A]{fontenc}
\usepackage[cp1251]{inputenc}
\usepackage[english,russian]{babel}      % russian as main

\usepackage{amsmath}    %Not needed with amsart
\usepackage{amsthm}     %Not needed with amsart
\usepackage{amsfonts,amssymb}

\usepackage{fullpage}

\usepackage{caption}
\usepackage{xcolor}
\usepackage{graphicx}
\usepackage{wrapfig}

\makeatletter
\def\ps@pprintTitle{%
     \let\@oddhead\@empty
     \let\@evenhead\@empty
     \def\@oddfoot{}%{\footnotesize\itshape
     \let\@evenfoot\@oddfoot}
\makeatother

\DeclareMathOperator{\const}{const}

\begin{document}

\setlength{\abovedisplayskip}{6pt plus 2pt minus 3pt}  % 8 - 2 - 4
\setlength{\belowdisplayskip}{6pt plus 2pt minus 3pt}  %8 - 2 - 4

\def\refname{Литература}%{Список литературы}%
\def\bibname{Литература}%{Библиография}%

\newsavebox{\TmpBox}
\newcommand{\tmp}{}
\newlength{\tmplength}                    %
\newcommand{\Mark}[1]{{\small\{\underline{#1}\}}\marginpar{??}}

\newtheorem*{thm}{Теорема}
\newtheorem{prop}{Утверждение}

%
%        Labelled/unlabelled equations
%
\newcommand{\Equa}[2]{\begin{equation}#2\label{#1}\end{equation}}
%        NonLabelled equation
\newcommand{\equa}[1]{\[ #1 \]}

%\newcommand{\eqref}[1]{{\rm(\ref{#1})}}
%\newcommand{\eqrefeq}[2]{{\rm(\ref{#1},\ref{#2})}}
%
%      Fractions
%
\newcommand{\myfrac}[2]{{\ifmmode{}^{#1}\!/_{\!#2}\else${}^{#1}\!/_{\!#2}$\fi}}
%
%        Math. functions
%
\newcommand{\abs}[1]{\left\lvert#1\right\rvert}
\renewcommand{\Re}[1]{\mathop{\rm Re}\nolimits\,(#1)}
\renewcommand{\Im}[1]{\mathop{\rm Im}\nolimits\,(#1)}
\newcommand{\sign}{\mathop{\rm sgn}\nolimits}
\renewcommand{\tan}{\mathop{\rm tg}\nolimits}
\renewcommand{\cot}{\mathop{\rm ctg}\nolimits}
\renewcommand{\arctan}{\mathop{\rm arctg}\nolimits}
\newcommand{\atanh}{\mathop{\rm arth}\nolimits}
\renewcommand{\cosh}{\mathop{\rm ch}\nolimits}
\newcommand{\acosh}{\mathop{\rm Arch}\nolimits}
\newcommand{\asinh}{\mathop{\rm Arsh}\nolimits}
\renewcommand{\sinh}{\mathop{\rm sh}\nolimits}
\renewcommand{\tanh}{\mathop{\rm th}\nolimits}
\newcommand{\Deg}[1]{{\ifmmode{#1}^\circ\else${#1}^\circ$\fi}}
\newcommand{\e}{{\mathrm{e}}}
\newcommand{\Exp}[1]{\e^{#1}}
%

%       sqrt(-1):
\newcommand{\iu}{{\mathrm{i}}}
%
%        Derivatives, integral
\newcommand{\diff}[1]{{\mathrm d}#1}
\newcommand{\D}[2]{#1_{#2}^{\prime}}
\newcommand{\DD}[2]{#1_{#2}^{\prime\prime}}
\newcommand{\Int}[4]{\displaystyle\int\limits_{#1}^{#2}{#3}\,{\mathrm d}#4}
\newcommand{\IntF}[5]{\displaystyle\int\limits_{#1}^{#2}\frac{#3\,{\mathrm d}#5}{#4}}
\newcommand{\Dfrac}[2]{\dfrac{{\mathrm d}#1}{{\rm d}#2}}
\newcommand{\Pd}[2]{\frac{\partial#1}{\partial#2}}
\newcommand{\ieq}{\,{=}\,}
\newcommand{\ineq}{\,{\not=}\,}
\newcommand{\ilt}{\,{<}\,}
\newcommand{\igt}{\,{>}\,}
\newcommand{\In}{\,{\in}\,}
\renewcommand{\le}{\leqslant}
\newcommand{\ile}{\,{\le}\,}
\renewcommand{\ge}{\geqslant}
\newcommand{\ige}{\,{\ge}\,}

\newcommand{\Eqref}[1]{\stackrel{\eqref{#1}}{=}}
\newcommand{\eqdef}{\stackrel{def}{=}}
%
%  > 0, < 0, >= 0,  etc   for in-line formulae: no $...$-delimiters required !
%    \GT{F(x,y)}  ==  $ F(x,y)\,{>}\,0 $
%
%
\newcommand{\GT}[1]{$#1\igt0$}
\newcommand{\GE}[1]{$#1\ige0$}
\newcommand{\LT}[1]{$#1\ilt0$}
\newcommand{\LE}[1]{$#1\ile0$}
\newcommand{\EQ}[1]{$#1\ieq0$}
\newcommand{\NE}[1]{$#1\ineq0$}
%
%            Others
\newcommand{\So}{\quad\Longrightarrow\quad}
\renewcommand{\Vec}[1]{\overrightarrow{#1}}
\newcommand{\HM}{\hphantom{{-}}}

%
%  Refs to figs.
%
\def\FigDir{}

\newcommand{\Figref}[1]{\ref{F#1}}
\newcommand{\Reffig}[2][]{рис.~\Figref{#2}\textit{#1}} 
\newcommand{\RefFig}[2][]{Рис.~\Figref{#2}\textit{#1}} 
\newcommand{\RefFigE}[2][]{Fig.\,\Figref{#2}\textit{#1}} 

\newcommand{\Infigw}[2]{%  Width Id
\includegraphics[width=#1]{\FigDir#2.eps}}

\newcommand{\Infig}[3]{\Infigw{#1}{#2}\caption{#3}\label{F#2}}

%                                        +-------+ 
%   +------...-----+   +------...-----+  |       | 
%   |              |   |              |  |       |
%   +------...-----+   +------...-----+  +--...--+
%       Fig. 1.           Fig. 2.          Fig. 3.

%
%          +------...-----+
%          |              |
% Fig. 1.  +------...-----+
%
%
\newsavebox{\CapBox}
\newcommand{\Lwfig}[3]{%  {t,b,h} {space_after_pic/-1pt}  {id} (no caption)
\begin{figure}[#1]
\settowidth{\tmplength}{Рис.~9.99.}%
\savebox{\CapBox}{\parbox{\tmplength}{\caption{}\label{F#3}}}%
\addtolength{\tmplength}{3mm}%
\raisebox{5pt}{\parbox[b]{\tmplength}{\usebox{\CapBox}}}%
%  Special meaning of -1pt: try to center picture
\ifdim#2=-1pt
  \addtolength{\tmplength}{\tmplength}%
\fi%
  \addtolength{\tmplength}{-\textwidth}%
  \setlength{\tmplength}{-\tmplength}%
%\else
  \addtolength{\tmplength}{-#2}%
%\fi
\includegraphics[width=\tmplength]{\FigDir#3.eps}
%\vspace{3mm}
\end{figure}
}

\newcommand{\Cwfig}[4]{% {h}  {width}  {id} {caption}
\begin{figure}[#1]
\centerline{\includegraphics[width=#2]{\FigDir#3.eps}}\caption{#4}\label{F#3}
\end{figure}
}

\newcommand{\Cfig}[4]{%  [pos] {width}  {id}  (caption)
\begin{figure}[#1]%
%\captionstyle{centerlast}
\centering%
\Infig{#2}{#3}{#4}%
%Infig{\ifdim#2=0pt\textwidth\else#2\fi\relax}{#3}{(#3) #4}
\end{figure}}

%     +------------+             
%     |            |             
%     +------------+             
%      Fig. 1. Text                
\newcommand{\Bfig}[3]{%  {width}  {id}  (caption)
\parbox[b]{#1}{\Infig{#1}{#2}{#3}}%\HideDisplacementBoxes
}

\newcommand{\Ffig}[3]{%  "here" Id  Caption
   \begin{figure}[#1]%
      \Infig{\textwidth}{#2}{#3}%
   \end{figure}%
}

% For bibliography:
%
\newcommand{\Jnum}[1]{$N$#1}
\newcommand{\Jitem}[6]{%    {author.}{title}{journal}{vol}{year}{N, pages}
#1 #2. {\it #3,\/} {\bf#4}(#5), #6}
%
%------------------------------------------- end of MyMath.tex

\newcommand{\quo}[1]{<<#1>>}
%\newcommand{\quo}[1]{\glqq#1\grqq}            % internal
%\Rus 

%+Title
\author{Alexey Kurnosenko}
\title{On tractrices of planar curves}
\date{}%{\today}
%\email{AlexeyKurnosenko@gmail.com}

%\English

\begin{abstract}
{\color{blue}\em The article is written in Russian.}\\
Tractrices of planar curves, in particular, a family of
tractrices of a circle, are considered.
Some new observations (including arc-length parametrization, Chezaro equation) 
and corrected reference informations are provided.
\end{abstract}

\maketitle\thispagestyle{plain}

These notes have been written after a rather unsuccessful search for detailed information 
on tractrices of a circle. 
The most comprehensive description was found in the well-known monograph of Gino Loria~\cite{Loria}.
We derive the natural parametrization of these curves and general Chezaro equation.
``Internal tractix'', omitted in~\cite{Loria}, is added to the classification of tractrices,
proposed herein.
%Some misprints from~\cite{Loria}, reproduced in modern handbooks, are corrected.

Let the point $A$ move along a directed circle (whose curvature~$K\lesseqgtr0$
defines the direction of movement), 
pulling the mass point~$B$, attached to $A$ by the thread of the length~$T$.
Tractrix is the trajectory of~$B$. 
The Chezaro equation of tractrix of the circle is
\equa{%TrxKs}{%
      k^\pm(s;w,T)=\pm\frac{w-(1{+}w)\Exp{-\frac{s}{T}}}%
                    {T\sqrt{1-\left[w-(1{+}w)\Exp{-\frac{s}{T}}\right]^2}},
            \quad \begin{array}{ll}where &w=\pm KT,\\and& w > {-1}\end{array}
}
(plotted in \RefFigE{TrxKs}). The following cases are described ($R\ieq\abs{K}^{-1}$):
\renewcommand{\tmp}[4]{\mbox{\bf#1:}&\mbox{#2,}&#3 & \text{#4}}%
\equa{%
\begin{array}{llrr}
\tmp{T1}{ external tractrix with a long thread $T>R$}{w>1\HM{}}{(\RefFigE[a]{TrxRw});}\\
\tmp{T2}{ external tractrix with $T\ieq R$ (polar, spiral tractrix)}{w=1\HM{}}{(\RefFigE[b]{TrxRw});}\\
\tmp{T3}{ external tractrix with $T<R$}{0< w< 1\HM{}}{(\RefFigE[c]{TrxRw});}\\
\tmp{T4}{ tractrix of a straight line ($R\ieq\infty$)}{w=0\HM{}}{(\RefFigE[d]{TrxRw});}\\
\tmp{T5}{ internal tractrix ($T<R$)}{-1< w< 0\HM{}}{(\RefFigE[e]{TrxRw});} %\\
%\tmp{T6}{ reverse tractrix, $T<{-R}<0$}{w<{-1}}{(\RefFigE[г]{TrxRw}).}
\end{array}
}
(the case  $w<{-1}$, $T<{-R}<0$ produces reverse tractrix, like $GF$ in \RefFigE{TrxRev}).
There is also a trivial case (\RefFigE[d]{TrxDef}) when tractrix is a circle (or straight line).
In \RefFigE{TrxOrtho} two tractrices of the curve $A_1A_2A_3$ are shown as the otrhogonal trajectories
of the family of circles, centered along $A_1A_2A_3$.

\RefFigE{TrxArch} shows involute $E_1E_2E_3$ of the circle of radius $R$
as a tractrix of Archimedian spiral $A_1A_2A_3$.
\RefFigE{TrxPer} illustrates periodic tractrices of periodic curves.
Polar equations of tractrices are given in \eqref{Polar} 
(two misprints from~\cite{Loria} in polar equations $[p(t),\,\varphi(t)]$ 
are underbraced in Sec.\,\ref{sec:err}).

\renewenvironment{abstract}{\global\setbox\absbox=\vbox\bgroup
  \hsize=\textwidth\def\baselinestretch{1}%
  \noindent\unskip\textbf{Аннотация}
 \par\medskip\noindent\unskip\ignorespaces}
 {\egroup}

\Russian

%\makeatletter
%\let\@author{А.\,И.~Курносенко}  %\@empty
\def\elsauthors{}
%\makeatother
%+Title
\author{Алексей Курносенко}
%\title{О трактрисах плоских кривых}
\title{Трактат о трактрисах}
\date{}%{\today}
%\ead{AlexeyKurnosenko@gmail.com}
%-Title

%??  УДК~514.752.22

%\begin{flushright}
%УДК 514.752.22   % (Кривые на плоскости) 
%\end{flushright}

%+Abstract
\begin{abstract}
Эти заметки посвящены трактрисам плоских кривых, в частности~---  семейству 
трактрис окружности. 
Приводятся новые наблюдения (периодические трактрисы, натуральное уравнение трактрис окружности),
справочная информация об этих кривых, предложена их классификация.
Указаны неточности, перекочевавшие из классического описания~\cite{Loria} в современные справочники.
Написано в духе курсовой работы старательного студента.
\end{abstract}

%-Abstract

\maketitle\thispagestyle{plain}

Если один конец~($A$) нерастяжимой нити $AB$ %длины $T$ 
вести вдоль кривой~$Z$,
то тяжёлая материальная точка, прикреплённая к концу~$B$,
опишет {\em трактрису} кривой~$Z$ (\Reffig{TrxDef}).
Кривая~$Z$ по отношению к трактрисе является 
{\em эквитангенциальной кривой}.

Часто просто трактрисой называют называют трактрису прямой линии 
(\Reffig[b]{TrxDef}). Эта кривая хорошо известна и подробно описана.
Описания трактрис окружности более скупы и местами ошибочны.
Обычно рассматривается трактриса c длиной поводка,
равной радиусу окружности, называемая {\em спиральной} или 
{\em полярной} трактрисой (\Reffig[c]{TrxDef}). Для этих двух кривых
известны и натуральные уравнения, частные случаи полученного ниже общего уравнения~\eqref{TrxKs}.

\Cfig{b}{.88\textwidth}{TrxDef}{Образование трактрисы кривой $Z$ (окружности, прямой)}%

Видимо, наиболее полное описание трактрис окружности дано
в классическом труде~\cite{Loria}, однако в нём допущены неточности,
перекочевавшие и в другие справочники.
Пытаясь приостановить их распространение, 
мы обратим на них особое внимание (п.\,\ref{sec:err}).
История исследования этих кривых, свойства и приложения
трактрисы прямой описаны в~\cite{Savelov}.\smallskip

%\Cfig{htb}{.8\textwidth}{TrxSolo}{TrxSolo}%
\Cfig{t}{.9\textwidth}{TrxEBTL}{Трактрисы некоторых кривых (лемниската, восьмёрка, синусоида,
составной контур, пятиугольник, квадрат). Чёрная кривая --- ведущая (её ориентация указана).}%

На \Reffig{TrxEBTL} показаны трактрисы нескольких кривых.
Сплошные векторы указывают начальное положение ведомой точки, пунктирные --- положение в точке останова
(если произошло ослабление нити).

\section{Трактриса произвольной кривой}

\subsection{Построение трактрисы}
%\Wfigr{8}{.5\textwidth}{TrxKin}{}%
Пусть $u(l),v(l)$~--- текущее положение ведущей точки~$A$ на данной кривой, $l$~--- натуральный параметр.
Угол $\theta(l)=\arg(u'_l+\iu v'_l)$ определяет направление её касательной (\Reffig{TrxKin}). 
Точки искомой трактрисы обозначим $x(l),\,y(l)$; 
\ $\tau(l)=\arg(x'_l+\iu y'_l)$~--- направление касательной к ней,
совпадающее с направлением натянутой нити (вектора $\Vec{BA}$).
Тогда при длине нити~$T$
\Equa{TrxAny}{%
%   \left\{\!  
   \begin{array}{l}
      x(l)=u(l)-T\cos\tau(l),\\
      y(l)=v(l)-T\sin\tau(l).
   \end{array}   
}
Натуральный параметр трактрисы обозначим~$s$.
Дифференцируем по~$l$:
\equa{
   \begin{array}{l}
      \Dfrac{x}{s}\Dfrac{s}{l}=\Dfrac{u}{l}+T\sin\tau\Dfrac{\tau}{l},\\[8pt]
      \Dfrac{y}{s}\Dfrac{s}{l}=\Dfrac{v}{l}-T\cos\tau\Dfrac{\tau}{l},
    \end{array}
    \quad \Longrightarrow \quad
   \begin{array}{l}
      \cos\tau\Dfrac{s}{l}=\cos\theta+T\sin\tau\Dfrac{\tau}{l},\\[8pt]
      \,\sin\tau\Dfrac{s}{l}=
                               \sin\theta -T\cos\tau\Dfrac{\tau}{l}.
    \end{array}
}
\begin{figure}[t]
\Bfig{0.48\textwidth}{TrxKin}{К выводу натурального уравнения трактрисы}%
\hfill
\Bfig{0.47\textwidth}{TrxArch}{Эвольвента окружности  $(E_1E_0E_2)$ как трактриса 
спирали Архимеда $(A_1A_0A_2)$}%
\end{figure}
%
%
%\noindent%
%\Wfigr{9}{.5\textwidth}{TrxKin}{}%
Решая линейную относительно $\D{\tau}{l}$, $\D{s}{l}$ систему, получаем
\begin{eqnarray}
        \D{\tau}{l}&=&\frac{\sin(\theta{-}\tau)}{T},\label{TrxAny1}\\
        \D{s}{l}&=&\cos(\theta{-}\tau)\label{TrxAny2}.
\end{eqnarray}
Решения~$\tau(l)$ уравнения~\eqref{TrxAny1} с известной функцией~$\theta(l)$
%и начальным условием $\tau(0)\ieq\tau_0$ 
достаточно для построения трактрисы по формулам~\eqref{TrxAny}.
Проинтегрировав затем~\eqref{TrxAny2}, мы находим $s(l)$ 
и обратную функцию~$l(s)$, что позволяет перейти 
к натуральной параметризации трактрисы.

Вошедший в эти уравнения угол между направлениями 
движения двух точек, $\theta{-}\tau$,
должен быть острым, иначе нить не будет тянуть точку~$B$: 
\Equa{Defnu}{%
         \nu(l)=\theta(l){-}\tau(l),\qquad 
         -\frac{\pi}{2}\ile\nu(l)\ile\frac{\pi}{2}. 
}
При $\nu\ieq{\pm\pi/2}$ траектория начинается или прекращается.
Подстановка~\eqref{Defnu} оставляет в уравнении~\eqref{TrxAny1}
переменные, не зависящие от выбора системы координат: 
\Equa{TrxAnyDE}{%
    \Dfrac{\nu}{l}=q(l)-\frac{\sin\nu(l)}{T},
}
где $q(l)\ieq\D{\theta}{l}(l)$ --- кривизна кривой $[u(l),v(l)]$.
Натуральное уравнение искомой трактрисы получим так:
\Equa{TrxAny3}{%
  k(s)=\Dfrac{\tau}{s} =\Dfrac{\tau}{l}{\cdot}\Dfrac{l}{s}
     % =\frac{1}{T}\frac{\sin(\theta{-}\tau)}{\cos(\theta{-}\tau)}
      =\frac{\sin(\theta{-}\tau)}{T}{\cdot}\frac{1}{\cos(\theta{-}\tau)}
      =\frac{1}{T}\tan\nu(l(s)).
}

Существование, единственность и непрерывная зависимость решений уравнения~\eqref{TrxAnyDE}
от начальных условий обеспечивается выполнением условия Липшица:
\Equa{Lipshiz}{%
    %\abs{\frac{T q(l){-}\sin\nu_1}{T}-\frac{T q(l){-}\sin\nu_2}{T}}=
     \abs{\Dfrac{\nu_1}{l}-\Dfrac{\nu_2}{l}}=
     \frac{2}{T}\abs{\sin\frac{\nu_1{-}\nu_2}{2}\cos\frac{\nu_1{+}\nu_2}{2}}
     \le \frac{\abs{\nu_1{-}\nu_2}}{T}\,.
}
\begin{figure}[t]
\Bfig{.66\textwidth}{TrxPush2}{Трактрисы типа \quo{тяни-толкай} для эллипса и квадрата}%
\hfill%
\Bfig{.25\textwidth}{TrxRev}{Звёздчатые тракт\-рисы окружности}%
\end{figure}

\subsection{Трактрисы типа \quo{тяни-толкай}}

%\noindent{\bf Трактрисы типа ``тяни-толкай''.\ }
Теперь заменим нить стержнем,
полагая, что при достижении положения, в котором нить бы ослабла, 
стержень будет продолжать работать в режиме толкания ведомой точки.
Будем считать, что толкание происходит по траектории,
которая образовалась бы при выполнении обратного движения в режиме \quo{тяни},
и рассматривать траекторию толкания как реверсную по отношению к трактрисе с нитью.
Тогда мы остаёмся в рамках модели~\eqref{TrxAny}
с заменой~$T$ на кусочно-постоянную функцию ${T}\sign\cos\nu(l) \ieq{\pm T}$,
отрицательную в режиме толкания.

Кривая $A_1A_2A_3$ на \Reffig{TrxArch}~--- спираль Архимеда $p(\varphi)=R\varphi$
(ветвь $A_1A_2$ соответствует \LT{\varphi}).
Легко убедиться, что при длине нити-стержня $T\ieq R$ получим в качестве одной из трактрис 
кривую $E_1E_2E_3$~--- эвольвенту окружности радиуса $R\ieq \abs{A_2E_2}$.
Режим \quo{тяни} показан на рисунке пунктирными стрелками ($E_1A_1$);
их направление совпадает с касательной к трактрисе. Переключение режима
происходит в точке $E_2/A_2$, бывшей точке останова, теперь~--- точке возврата. 
Направление касательной меняется в ней скачком.

Примеры трактрис типа \quo{тяни-толкай} приведены также на \Reffig{TrxPush2} и \Figref{TrxRev}.

\subsection{Теорема о периодических трактрисах} 
Следующее утверждение мы докажем в ослабленном варианте, 
ограничив кривизну ведущей кривой, а проиллюстрируем без этого ограничения (\Reffig{TrxPer}). 

\begin{thm}
У периодической кривой ограниченной кривизны существует семейство
периодических трактрис.
\end{thm}

%\noindent {\bf Доказательство.\ } 
\begin{proof}[\bf\proofname]
Если кривизна ограничена, т.е. $\abs{q(l)}\ilt\infty$,
то длину нити можно ограничить минимальным радиусом кривизны:
\Equa{TrxkT}{%
       T\ilt \min\abs{1/q(l)},\text{~~~т.\,е.~~~} \abs{q(l)\,T}<1.
}
Тогда, если $\nu(l)$ возрастает, приближаясь слева к~$+{\pi}/{2}$, 
то \LT{\D{\nu}{l}},  и $\nu(l)$~начнёт убывать~\eqref{TrxAnyDE}. 
Приближаясь справа к~$-{\pi}/{2}$, убывающая функция $\nu(l)$ становится возрастающей.
Таким образом, при начальных значениях $\nu(0)\In[-{\pi}/{2},{\pi}/{2}]$
решения  $\nu(l)$ 
более не выходят за пределы интервала $(-{\pi}/{2},{\pi}/{2})$, 
и мы всё время остаёмся в режиме тянущей нити. 

\Cfig{t}{.9\textwidth}{TrxPer}%
{Семейства периодических трактрис для четырёх периодических кривых}%

Пусть теперь $\nu_1(l)$ и $\nu_2(l)$~--- два решения 
уравнения~\eqref{TrxAnyDE} 
с начальными условиями $\nu_1(0)<\nu_2(0)$, \ $\nu_{1,2}(0)\in[-\pi/2,\pi/2]$. 
Тогда для функции $\delta(l)=\nu_2(l){-}\nu_1(l)$ справедливо
\equa{
     \delta(0)>0\text{~~~и~~~}
     \D{\delta}{l}=-\frac{1}{T}[\sin\nu_2(l){-}\sin\nu_1(l)]=
                   -\frac{2}{T}\sin\frac{\delta(l)}{2}\cos\frac{\nu_2(l){+}\nu_1(l)}{2}<0.
}
Ситуация ${\exists l_0}:\:\delta(l_0)\ieq 0$ невозможна, так как она бы означала
нарушение единственности решения в окрестности точки $l\ieq l_0$~\eqref{Lipshiz}.
Следовательно, \GT{\delta(l)}, 
%производная $\delta'(l)$ отрицательна,
и $\delta(l)$ монотонно убывает.
Два решения  $\nu_1(l)$ и $\nu_2(l)$ сближаются, оставаясь внутри
интервала $\left(-{\pi}/{2},{\pi}/{2}\right)$.
При $\nu_1(0)\ieq{-\pi/2}$ и $\nu_2(0)\ieq{\pi/2}$ мы получаем
непрерывное и сжатое отображение отрезка $\left[-{\pi}/{2},{\pi}/{2}\right]$  на 
отрезок $\left[\nu_1(L),\nu_2(L)\right]$,
где~$L$~--- период функции~$q(l)$.
Оно имеет неподвижную точку~$\nu^\star$.
Взяв её в качестве начального условия,
получим искомое периодическое решение уравнения~\eqref{TrxAnyDE}. 
Параметром полученного семейства служит длина нити~$T$.%
\end{proof}

%\begin{figure}[!t]
%\Hfig{0mm}{TrxPer}%
%\\~~\\~~\\%
%\Hfig{0mm}{Trx90}%
%\end{figure}

Условие ограниченности кривизны \eqref{TrxkT}  на~\Reffig{TrxPer} соблюдено не везде.
Кривизны эллипса и восьмёрки ограничены, но в изображённые семейства
включены также периодические трактрисы со значениями~$T$, нарушающими условие~\eqref{TrxkT}.
Гладкая кривая $E_0E_1E_2$ включает полуокружность $E_1E_2$ и две дуги
эвольвенты окружности, $E_0E_1$  и $E_2E_0$, 
\ с \ $k(E_0)\ieq\infty$. У ведущей кривой $F_0F_1F_2$ имеются изломы,
которые можно представить как импульсы кривизны (в виде дельта-функций Дирака).
%И то, и другое не вписывается в~\eqref{TrxkT}.
В иллюстрациях присутствуют режимы толкания, не рассмотренные в доказательстве.

Периодические трактрисы окружности представлены как тривиальная
трактриса на \Reffig[d]{TrxDef} и звёздчатые кривые на \Reffig{TrxRev}.

\section{Трактрисы окружности и прямой}

Джино Лориа выводит дифференциальное уравнение трактрисы окружности радиуса~$a$
с поводком длины~$l$ в полярных координатах $(\varrho,\omega)$ в виде
\equa{%
   \varrho^2+\frac{2l\varrho\,\diff{\varrho}}%
   {\sqrt{(\diff{\varrho})^2+\varrho^2(\diff{\omega})^2}}=a^2-l^2.      
}
Мы получим полярные и параметрические уравнения этих кривых исходя из натурального уравнения
\Equa{TrxKs}{%
      k^\pm(s;w,T)=\pm\frac{w-(1{+}w)\Exp{-\frac{s}{T}}}%
                    {T\sqrt{1-\left[w-(1{+}w)\Exp{-\frac{s}{T}}\right]^2}},\quad
                    \begin{array}{ll}\text{где} & w={\pm KT},\\ \text{и}& w>{-1}.\end{array}
}
Здесь \GT{T} --- длина поводка, $K\lesseqgtr 0$~--- кривизна окружности,
знак которой определяет ориентацию окружности
(и тем самым направление движения ведущей точки). 
При этом, если поводок короче радиуса окружности, $T\ile\abs{K}^{-1}$, мы можем, в зависимости от
начальных условий, получить трактрисы двух типов, ниже названные внешней и внутренней 
(есть и третья возможность~--- ``тривиальная'' трактриса, \Reffig[d]{TrxDef}).
При $T\ige\abs{K}^{-1}$ только один из вариантов $w\ieq{\pm KT}$ удовлетворяет условию 
 $w>{-1}$, и существует единственная (внешняя) трактриса.

\subsection{Вывод натурального уравнения}

Пусть ведущая кривая --- окружность
кривизны $q(l)\ieq K\ieq \const$ (радиуса $R\ieq\abs{K}^{-1}$), 
что включает и случай прямой, \EQ{K}.
Пeрейдём в~\eqref{TrxAnyDE} от длины дуги~$l$ окружности к 
длине дуги~$s$ трактрисы:
\Equa{TrxDEs}{%
    \D{\nu}{s}=\D{\nu}{l}\cdot\D{l}{s}=\left(K-\frac{\sin\nu}{T}\right)\cdot\frac{1}{\cos\nu},
    \qquad\mbox{т.е.}\qquad
    \Dfrac{\nu}{s}=\frac{KT-\sin\nu(s)}{T\cos\nu(s)}.
}
Уравнение~\eqref{TrxDEs} при $\abs{KT}\ilt 1$ имеет тривиальное решение 
\Equa{TrxTriv}{%
   KT{-}\sin\nu(s)\equiv 0\qquad\Longrightarrow\quad
    \begin{array}{l}
       \nu(s)= \arcsin KT = \nu^\star = \const ,\\
       k(s) \Eqref{TrxAny3}\frac{1}{T}\tan\nu^\star= \const,
    \end{array}
}
означающее движение ведомой точки по прямой (если \EQ{K}) либо по окружности радиуса 
$\sqrt{R^2{-}T^2}$ (\Reffig[d]{TrxDef}).

Перемена знака выражения $KT{-}\sin\nu(s)$ в некоторой точке
% $s\ieq s_0$
невозможна, ибо привела бы к переходу 
на тривиальную интегральную кривую~\eqref{TrxTriv}
с нарушением единственности решения 
в окрестности этой точки. Из знакопостоянства числителя 
и всей правой части~\eqref{TrxDEs}
следуют монотонность функции~$\nu(s)$ и общее решение
\equa{%TrxInt}{%
      s=\Int{\nu(0)}{\nu(s)}{\frac{T\cos\nu}{KT-\sin{\nu}}}{\nu}
       =-T\ln\frac{KT-\sin{\nu(s)}}{KT-\sin{\nu(0)}}\So
       \sin\nu(s) = KT-\left[KT-\sin\nu(0)\right]\Exp{-\frac{s}{T}}.
}
%\smallskip

{\bf 1.}~Пусть $\nu(s)$~возрастает. 
%  \GT{T}, \GT{KT{-}\sin\nu(s)};
Тогда, если $\nu(0)\igt{-\pi/2}$, 
можно положить $\nu(0)\ieq{-\pi/2}$, 
включив тем самым предысторию движения. Обозначив $w = {+}KT$, 
%из \eqref{TrxInt} 
получаем 
\Equa{Nus}{%
        \sin\nu(s)=w-(1{+}w)\Exp{-\frac{s}{T}},\qquad 
        \cos\nu(s)=\sqrt{1-\sin^2\nu(s)}.
}
Натуральное уравнение \eqref{TrxAny3} для трактрисы окружности принимает вид $k^+(s;w,T)$
из~\eqref{TrxKs}.
%\equa{%TrxKs}{%
%      k(s;w,T)=\frac{w-(1{+}w)\Exp{-\frac{s}{T}}}%
%                    {T\sqrt{1-\left[w-(1{+}w)\Exp{-\frac{s}{T}}\right]^2}}.
%            %\quad w={\pm\frac{T}{R}},\quad w\ineq {-1}
%}
Отметим, что для \GE{s} радикал в знаменателе действителен при
\Equa{TrxwT}{%
    \begin{array}{rcl}
              w\igt 1,& \quad T\igt 0,\quad& 0\ile s\ile S_1,\\
       -1\ilt w\ile 1,& T\igt 0,& 0\ile s\ile \infty,\\
              w\ilt {-1},& T\ilt 0,& 0\ile s\ile S_1,
    \end{array}
    \qquad  \mbox{где~~~}S_1\ieq T\ln\frac{w+1}{w-1}.
}
Значение \LT{T} в третьей строке соответствует в исходных уравнениях~\eqref{TrxAny}
режиму толкания и трактрисе, реверсной по отношению к трактрисе из первой строки.%
\footnote{% 
Если имеется кривая $k(s),\;s\in[0,S_1]$, то реверсную кривую можно получить как
$-k(-s)$ либо, при конечной длине $S_1$, как $-k(S_1{-}s)$.
Легко убедиться, что
\equa{
      -k^+(S_1{-}s;w,T)\equiv k^+(s;{-w},{-T}).
}
}
\smallskip

\Cfig{t}{.66\textwidth}{TrxKs}{Зависимость $k^+(s;w,T)$}

{\bf 2.}~Случай убывающей функции~$\nu(s)$ сводится к предыдущему с помощью симметрии.
При симметрии меняют знак ориентированный угол $\nu(s)$ $[\nu(0)={+\pi/2}]$ 
и кривизна~$K$ ведущей окружности. 
Приняв в этой ситуации $w= {-KT}$, мы сохраним значение параметра формы~$w$
и получим, как и следовало ожидать при симметрии, вариант $k^-(s;w,T)$
с убывающей кривизной.

\Cfig{t}{\textwidth}{TrxRw}{Различные типы трактрис окружности}

\subsection{Классификация  трактрис окружности}

Иллюстрации (\RefFig{TrxKs}, \Figref{TrxRw}) соответствуют варианту $k^+$ в \eqref{TrxKs}
(возрастание кривизны и угла $\nu$).

Направим ведущую окружность (прямую) из начала координат по оси ординат (\Reffig{TrxRw}). 
Тогда $\theta(0)\ieq\pi/2$, и $\theta(l)\ieq\pi/2+K l$.
Начальному условию $\nu(0)\ieq{-\frac{\pi}2}$ соответствуют 
направление касательной к трактрисе $\tau_0\Eqref{Defnu}\theta(0){-}\nu(0)\ieq\pi$
и координаты ведомой точки $(T,0)$.

Окружности с положительной кривизной \GT{K} уходят в левую полуплоскость (\Reffig[a,b,c]{TrxRw}).
Длина нити при \GE{K} может быть произвольной: \GE{w}. 

При \LT{K} окружности уходят в правую полуплоскость (\Reffig[e]{TrxRw}).
Ведомая точка должна находиться внутри окружности, 
а нить должна быть короче радиуса окружности 
$T\ilt R\ieq {-1/K}$ (иначе изначальное натяжение нити сразу ослабнет).
У трактрис этого типа \ ${-1}\ilt w \ilt 0$. 

Далее, игнорируя тривиальую трактрису~\eqref{TrxTriv}, 
мы рассматриваем пять разновидностей трактрис окружности/прямой:
%помеченных в тексте как T1--T5.
\renewcommand{\tmp}[3]{\mbox{\bf#1:}&\mbox{#2,}&#3}
\equa{%
\begin{array}{llr}
\tmp{T1}{ внешняя трактриса на длинном поводке, $T>R$}{w>1;\HM{}}\\
\tmp{T2}{ внешняя трактриса с $T\ieq R$ (полярная, спиральная)}{w=1;\HM{}}\\
\tmp{T3}{ внешняя трактриса на коротком поводке $T<R$}{0< w< 1;\HM{}}\\
\tmp{T4}{ трактриса прямой (поводок \quo{короткий}: $T\ilt R\ieq\infty$)}{w=0;\HM{}}\\
\tmp{T5}{ внутренняя трактриса (поводок короткий)}{-1< w< 0.\HM{}}%\\
%\tmp{T6}{ реверсная трактриса, $T<{-R}<0$}{w<{-1}.}
\end{array}
}

Длина трактрис T2--T5 неограничена.
%Трактриса T3 наматывается на окружность~\eqref{TrxTriv}.
%Трактриса T5 наматывается на неё изнутри.
Окружность (прямая)~\eqref{TrxTriv} служит асимптотической кривой для трактрис T3,\,T4,\,T5.
У трактрисы T2 асимптотическая окружноcть превращается в предельную  точку.

Три внешние трактрисы имеют точку перегиба:
\equa{%
    k(s_0;w,T)=0\quad\mbox{при}\quad s_0=T\ln\frac{w+1}{w}.
}
Касастельная к трактрисе в точке перегиба совпадает с касательной к окружности
в соответствующей точке ($\nu(s)=0\;\Rightarrow\; k(s)=0$).

\subsection{Параметрические уравнения}

Обозначения функций и параметров для конкретных типов мы иногда снабжаем индексом, 
например, $\tau_{1}(s)$, $w_5$. 

Далее мы используем замену переменных
\equa{%
  s\to \nu(s) \to \nu(s)+\frac{\pi}{2}=\xi(s)
}
и вводим тем самым будущие параметры кривых~$\xi$ и~$t$, 
\Equa{DefXi}{
      \begin{array}{rcl}
        \cos\xi(s)&\!\!\Eqref{Nus}\!\!&-\sin\nu= -w+(1{+}w)\Exp{-\frac{s}{T}},\\
        \sin\xi(s)&\!\!=\!\!&\HM{}\cos\nu=\sqrt{1-\cos^2\xi(s)},
      \end{array}\qquad
      t=\tan\frac{\xi(s)}{2}=  
    %\sqrt{\frac{(1+w)\left(1-\Exp{-s/T}\right)}{1{-}w+(1{+}w)\Exp{-s/T}}}\,,
     \sqrt{\frac{1-\Exp{-s/T}}{\frac{1{-}w}{1{+}w}+\Exp{-s/T}}}\,,
}
удобные тем, что в начале кривой \EQ{\xi(0)} и \EQ{t}.
Тогда 
\equa{%
    \diff{s}= \frac{T\sin\xi}{w{+}\cos\xi}\diff{\xi},\qquad
    k(s)=-\frac{1}{T}\ctg\xi(s).
}
Определим функцию $l(s)$, связывающую длины дуг двух кривых:
\renewcommand{\tmp}{\tan\frac{\xi(s)}{2}}
\equa{%
       l(s)\Eqref{TrxAny2}\IntF{0}{s}{}{\cos\nu(s)}{s}
       =\IntF{0}{s}{}{\sin\xi(s)}{s}
       =\IntF{0}{\xi(s)}{T}{w+\cos\xi}{\xi}
       =\frac{2T}{\sqrt{w^2{-}1}}%
          \arctan\left[\frac{w{-}1}{\sqrt{w^2{-}1}}\tmp\right].%
}
Полученное выражение действительно и при мнимых значениях радикалов ( $\abs{w}\ilt 1$) 
и, как предел, при $w\ieq 1$:
\Equa{DefPsi}{%           
   \begin{array}{l}
     l(s) = \left\{%
      \begin{array}{lclcl}
         l_{1}(s)&=&\frac{2T}{\sqrt{w^2{-}1}}%
            \arctan\left[\sqrt{\frac{w{-}1}{w{+}1}}\tmp\right],%
         &\quad& w>1;\\
          l_2(s)&=&T\tmp, 
         &&w=1;\\
         l_{3-5}(s)&=&\frac{2T}{\sqrt{1-w^2}}%
                 \atanh\left[\sqrt{\frac{1{-}w}{1{+}w}}\tmp\right],
          &&\abs{w}<1;  \\[6pt]
          l_{4}(s)&=& 2T\atanh\tmp\,,
%                 =T\ln\frac{1+\sin\xi(s)}{\cos\xi(s)}\,,
          &&w=0;
      \end{array}\right.\\ \\%[14pt]
      \psi(s)=K\,l(s) = \frac{w}{T}l(s).
    \end{array}  
}
В последней строке записан выраженный в угловых единицах путь, пройденный ведущей точкой по окружности.
При \GT{K} движение соответствует вращению 
полярного радиуса против часовой стрелки
(полюс в центре $(-R,0)$ ведущей окружности),
и $\psi_{1,2,3}(s)$ есть просто полярный угол точки~$A$
(\Reffig{TrxKin}, $\angle BOA$). При \EQ{K} имеем $\psi_4(s)\equiv 0$.
Отрицательное и убывающее значение $\psi_5(s)$ отражает 
убывание полярного угла точки~$A$, равного в этом случае $\pi{+}\psi(s)$
(полюс в центре $(R,0)$ ведущей окружности).
Теперь входящие в~\eqref{TrxAny} координаты $(u,v)$ точки~$A$
могут быть представлены как
\equa{%
   \begin{array}{rll}
      \mbox{T1,T2,T3:}&\left\{
         \begin{array}{lclcl}
          u(s)&=&R\cos\psi(s){-}R&=&\frac{T}{w}\left[\cos\psi(s){-}1\right],\\[6pt]
          v(s)&=&R\sin\psi(s)    &=&\frac{T}{w}\sin\psi(s)
         \end{array}\right.&
         \quad\left[R=\dfrac1K=\dfrac{T}{w}\right];\\[14pt]
      \mbox{T4:}&\left\{ 
         \begin{array}{lcl}
          u(s)&=&0,\\ 
          v(s)&=&l_4(s);
         \end{array}\right.&\\[12pt]
     \mbox{T5:}&\left\{
        \begin{array}{lclcl}
          u(s)&=&R\cos[\pi{+}\psi(s)]+R&=&\frac{T}{w}\left[\cos\psi(s){-}1\right],\\[6pt]
          v(s)&=& R\sin[\pi{+}\psi(s)]&=&\frac{T}{w}\sin\psi(s)
         \end{array}\right.&
         \quad\left[R=-\dfrac1K=-\dfrac{T}{w}\right];
  \end{array}
}
Направление $\tau(s)$ касательной к трактрисе равно
\equa{%
     \tau(s)=\tau(0)+\Int{0}{s}{k(\tilde s)}{\tilde s}
            =\pi+\Int{0}{\xi(s)}{\frac{-\cos\xi}{w{+}\cos\xi}}{\xi}
            %&=&\pi+\Int{0}{\xi(s)}{\left(\frac{w}{w{+}\cos\xi}-1\right)}{\xi}
            =\pi+\psi(s)-\xi(s).
}
Располагая выражениями для $u,v,\tau$, 
мы можем получить искомые функции $x(s),\:y(s)$  
непосредственно по фомулам~\eqref{TrxAny}
(т.е. без утомительного интегрирования $\cos\tau,\,\sin\tau$):
\begin{subequations}\label{TrxXY}
\begin{align}
   \mbox{T1,T2,T3,T5}{:}&
   \left\{
   \begin{array}{l}
      x(s)=T\left[\dfrac{\cos\psi(s){-}1}{w}+\cos[\psi(s){-}\xi(s)\right],\\[8pt]
      y(s)%&=&v(s)-\,T\sin\tau(s)\\
          =T\left[\dfrac{\sin\psi(s)}{w}+\sin[\psi(s){-}\xi(s)\right];
   \end{array}\right. \label{TrxXY1235}\\
   \mbox{T4}{:}&
   \left\{
   \begin{array}{l}
       x(s)=\HM{}T\cos\xi(s),\\
       y(s)=-T\sin\xi(s)+T\ln\dfrac{1+\sin\xi(s)}{\cos\xi(s)};
   \end{array}\right. \label{TrxXY4}
\end{align}
\end{subequations}
Эти параметризации используют выражения \eqref{DefXi} и \eqref{DefPsi} и соответствуют версии $k^+$ натурального 
уравнения~\eqref{TrxKs}.

\subsection{Трактриса прямой.}
Предельный переход $R{\to}\infty$ показан как  
\hbox{\em\RefFig[a]{TrxDef} ${\to}$ \RefFig[b]{TrxDef}}: 
внутренняя и внешняя трактрисы окружности становятся левой
и правой трактрисами ориентированной прямой.
Для правой трактрисы
\equa{
   k_4(s) = k^+(s;0,T)=\frac{-1}{T\sqrt{\Exp{\frac{2s}{T}}-1} },\qquad
    \begin{array}{l}
        \cos\xi(s)=\Exp{-\frac{s}{T}},\\
        \sin\xi(s)=\sqrt{1-\Exp{-\frac{2s}{T}} },\\
   \end{array}\qquad
    t=\sqrt{\frac{1-\Exp{-\frac{s}{T}}}{1+\Exp{-\frac{s}{T}}}},
}
а уравнения \eqref{TrxXY4} можно записать как
\equa{%
       x_4(s)=T\Exp{-\frac{s}{T}},\qquad
       y_4(s)=T\left(\atanh\sqrt{1-\Exp{-\frac{2s}{T}}} - 
                             \sqrt{1-\Exp{-\frac{2s}{T}}} \right),\qquad
       0\le s <\infty,         
}
либо
\equa{%
       x_4(t)=T\frac{1-t^2}{1+t^2},\qquad
       y_4(t)=T\ln\frac{1+t}{1-t}-\frac{2Tt}{1+t^2}\,,\qquad
       0\le t < 1.
}
Явное уравнение описывает и левую, и правую трактрисы (в положении на \Reffig[b]{TrxDef}):
\equa{%
       y(x)=T\left(\atanh\sqrt{1-\frac{x^2}{T^2}} - 
                   \sqrt{1-\frac{x^2}{T^2}} \right).
}
Чаще  трактриса прямой  рассматривается в таком положении:
\begin{picture}(45,14)(0,2)%
\put(0,0){\Infigw{45pt}{TrxLL}}
\end{picture}\,.

\Cfig{t}{.92\textwidth}{TrxMove}{Траектория полюса спиральной трактрисы при её качении по трактрисе прямой}

Упомянем здесь свойство, приведённое в~\cite{FrCite}: 
{\em если полярную трактрису с поводком длины $T_2$ 
катить по трактрисе некоторой прямой с поводком длины $T_4\ieq 2T_2$
совместив начальные положения,
то предельная точка полярной трактрисы
будет двигаться по этой прямой\,} (\Reffig{TrxMove}).

\subsection{Полярные уравнения}

Для вывода полярных уравнений перейдём в систему координат 
с началом в центре ведущей окружности (при этом в~\eqref{TrxXY1235} исчезнет слагаемое~$-1$): 
\equa{%
       p(\xi)=R\sqrt{1+2w\cos\xi+w^2},\quad R=\frac1{\abs{K}}=\abs{\frac{T}{w}}.
}
Здесь $p$~--- полярный радиус точки. Простые преобразования приводят к выражениям
(с заменой $\tan\frac{\xi}{2}\ieq t$):
\settowidth{\tmplength}{$0\ile t\ile\infty;$}
\Equa{Polar}{
%\begin{align}
\aligned
     p(t)&=          R\sqrt{1+2w\frac{1-t^2}{1+t^2}+w^2},&\\
     \varphi_{1}(t)&=
         \frac{2w}{\sqrt{w^2-1}}%
         \arctan\left(\sqrt{\frac{w{-}1}{w{+}1}}\,t\right)%
        -\arctan\left(\frac{w{-}1}{w{+}1}\,t\right)-\arctan t,\;
         &0\ile t\ile\infty;\\[6pt]
     \varphi_2(t)&=
          t-\arctan t, 
          %\qquad \left[p_2(t)\equiv p(t,1)=\frac{2T}{\sqrt{1+t^2}}\right],
          &0\ile t\ilt\infty;\\
     \varphi_{3,5}(t)&=
         \frac{2w}{\sqrt{1-w^2}}%
           \atanh\left(\sqrt{\frac{1{-}w}{1{+}w}}\,t\right)%
          + \arctan\left(\frac{1{-}w}{1{+}w}\,t\right)-\arctan t \equiv{}&\\
         &\equiv \frac{w}{\sqrt{1{-}w^2}}%
             \ln \frac{\sqrt{1{+}w} + t\sqrt{1{-}w}}{\sqrt{1{+}w} - t\sqrt{1{-}w}}
             +\arctan\frac{2t}{w{+}1 + (w{-}1)t^2},
             % &0\ile t\ilt \sqrt{\frac{1{+}w}{1{-}w}},\nonumber
            &\parbox{\tmplength}{$0\ile t\ilt \sqrt{\frac{1{+}w}{1{-}w}}\,.$}
\endaligned
}%\end{align}
Для~$\varphi_{5}(t)$ дополнительно введён поворот кривой на \Deg{180}  
в \quo{более каноничное} положение, такое, что $\varphi(0)\ieq 0$.
Трактриса T5 в представленном виде вращается
%, в отличие от остальных,
в сторону уменьшения полярного угла.

Положив для трактрисы T1 \ $w=\ch\omega,\;q=\th\frac{\omega}2$, получим
\equa{%
  p_1(t)=\frac{2R}{1-q^2}\sqrt{\frac{1+q^4t^2}{1+t^2}},\qquad
  \varphi_1(t)=\left(q+q^{-1}\right)\arctg(qt)-\arctg(q^2t)-\arctg t.
}
Положив для трактрис T3, T5 $w=\cos\omega,\;q=\tg\frac{\omega}2$, получим
\equa{%
  p_{3,5}(t)=\frac{2R}{1+q^2}\sqrt{\frac{1+q^4t^2}{1+t^2}},\qquad
  \varphi_{3,5}(t)=\left(q-q^{-1}\right)\atanh(qt)+\arctg(q^2t)-\arctg t.
}

Полярные уравнения можно представить в явном виде. 
Так, для спиральной трактрисы
\Equa{Trx2Pol}{%
    \varphi_2(p) = -\arctan\frac{\sqrt{4T^2-p^2}}{p} + \frac{\sqrt{4T^2-p^2}}{p}
           \equiv  -\arccos\frac{p}{2T}+\frac{\sqrt{4T^2-p^2}}{p}.
}
 
Трактриса T1 ограниченной длины~$S_1$~\eqref{TrxwT} ограничена и по~$\varphi$,
и заключена в полярном секторе ширины
\equa{%TrxDlt}{%
    \Delta\varphi=\varphi_1(\infty)-\varphi_1(0)=\varphi_1(\infty)
                 =\pi\left(\frac{w}{\sqrt{w^2-1}}-1\right).
}
На \Reffig{TrxRev} параметр~$w$ подобран так, чтобы образовались периодические звёздчатые кривые.

У трактрисы $GH$ на \Reffig{TrxRev} имеется симметричная ветвь $FG$ (\LT{t}),
являющаяся реверсной по отношению к трактрисе $GF$.
Вернувшись от параметра~$t\in[-\infty,+\infty]$ к $\xi\in[-\pi,\pi]$, 
кривую $FGH$ можно периодически продолжить на $\xi\in(-\infty,+\infty)$
и получить кривые в режиме \quo{тяни-толкай}.

\Cfig{t}{\textwidth}{TrxBeam}{Пучок окружностей, асимптотических линий пучка трактрис}%

\subsection{Пучок трактрис}

Семейство трактрис на коротком поводке $(\abs{w}\ilt 1)$
фиксированной длины~$T$ 
и одинаковыми начальными условиями $x(0),\,y(0),\, \tau(0)$
назовём {\em пучком трактрис\,} (\Reffig{TrxBeam}). 
Пунктиром на рисунке изображены асимптотические линии пучка~---
прямая и предельные окружности. 
% (на \Reffig{TrxRw} постоянным был параметр~$R$).
Оказывается, что
\smallskip\\%
%\noindent
{\em асимптотические линии пучка трактрис
образуют гиперболический пучок окружностей с расстоянием между фокусами $2T$. 
В одном из его фокусов находится общая начальная точка кривых семейства,
в другом~--- предельная точка трактрисы $T2$.
Асимптота трактрисы прямой служит радикальной осью пучка}.\smallskip

Действительно, у трактрисы T3 с параметром формы \GT{w} общий центр ведущей и предельной окружностей
расположен в точке $(-R,0)=(-T/w,0)$, у трактрисы T5 (\LT{w})~---
в точке $(R,0)=(-T/w,0)$. Радиусы предельных окружностей равны
$\sqrt{R^2{-}T^2}\ieq {T}\sqrt{1-w^2}/{\abs{w}}$.
Их неявные уравнения
\equa{%
       \left(x+\frac{T}{w}\right)^2+y^2-\frac{T^2(1-w^2)}{w^2}=0
}
преобразуются к виду $w(x^2+y^2+T^2)+2x=0$, 
включающему при \EQ{w} и ось ординат~--- асимптоту трактрисы прямой.
%У гиперболического пучка две общие точки мнимые.
%Легко убедиться, что для семейства \EQ{F(x,y;w)} 
%таковыми являются точки $(x_1,y_1)$ и $(x_2,y_2)$,
%где $x_{1,2}\ieq 0$, $y_{1,2}\ieq{\pm\iu T}$.

%любой пары \EQ{C(x,y;w')}, \EQ{C(x,y;w'')} окружностей семейства.
%\footnote{%
%Одно из определений радикальной оси --- прямая, проходящая через общие точки двух
%окружностей. В данном случае эти точки мнимые.
%}

\subsection{Ошибочные уравнения трактрис окружности}\label{sec:err}

Полярные уравнения трактрис приводятся в \cite{Loria} в виде
(с заменой обозначений $\theta\to t$, $\omega\to\varphi$, $\varrho\to p$)
\equa{%TrxLoria}{%
   \begin{array}{l}
      p(t)=\sqrt{a^2+l^2+2al\frac{1-t^2}{1+t^2}},\\[9pt]
     \varphi(t)=\left\{
     \begin{array}{ll}
        \arctan t - t &(a\ieq l),\\
        \underbrace{\arctan\frac{2lt}{(a{+}l)+(a{-}l)t^2}+\frac{2l}{n}}_{\mbox{\tiny ошибка!}}
               \arctan\left(\sqrt{\frac{l-a}{l+a}}\,t \right)\quad
               &(a\ilt l),\\[30pt]
        \arctan\frac{2lt}{(a{+}l)+(a{-}l)t^2}
               -\frac{l}{n}\ln\frac{a+l+nt}{a+l-nt}
               &(a\igt l),
     \end{array}
     \right.\\
     \mbox{где}\quad {n^2 = a^2{-}l^2}^{\strut}.
   \end{array}
}
%(при их цитировании в~\cite{Shikin} формула для~$n$ потерялась).
Это соответствует трактрисам T2, T1 и T3
с $a\ieq R$, $l\ieq T$.
% и убывающим полярным углом $\varphi(t)$.
Отметим, что в строке $a\ilt l$ значение~$n$ становится мнимым.
Но и с заменой~$n$ на \ $-\sqrt{l^2{-}a^2}$ остаётся ошибка
в первом слагаемом.
При переходе через 
%значение 
$t\ieq t_0\ieq\sqrt{\frac{l+a}{l-a}}$ оно претерпевает скачок от $+\pi/2$
до \ $-\pi/2$,
и на кривой образуется разрыв.\footnote{%
К ошибке привело, по-видимому, использование формулы
\equa{%AtanXY}{%
      \arctan x+\arctan y= \arctan\frac{x+y}{1-xy}\,,\quad
      \mbox{если}\;\abs{xy}\ilt 1\;\mbox{~либо~}\;\sign x\ineq \sign y,
} 
без учёта условий её применимости, нарушающихся при $t\ige t_0$.}

Здесь также упущена трактриса T5, которая получилась бы при дополнении
случая $a\igt l$ условиями \LT{a}, $\abs{a}\igt l$.

В \cite{Loria,Savelov,Shikin},
наряду с правильным полярным уравнением~\eqref{Trx2Pol} трактрисы T2
приводится ошибочный вариант: 
\equa{%TrxErr}{
    \varphi(p) = -\underbrace{\arcsin}_{\arccos\,!}%
    \frac{p}{2T}+\frac{\sqrt{4T^2-p^2}}{p}\qquad
    \left[\:\Longrightarrow\:
    \mbox{
     \begin{picture}(33,20)(0,12)%
     \put(0,0){\Infigw{33pt}{TrxShikin}}
     \end{picture}}
    \right].
}
%В уравнении следует заменить арксинус на арккосинус. 

\section{Трактрисы как ортогональные траектории}

На \Reffig{TrxOrtho} ведущая кривая $A_1A_2A_3$ составлена из параболы и прямой.
%(с нарушением гладкости в переходной точке). 
Показаны две трактрисы этой кривой.
Рисунок призван  проиллюстрировать следующий вполне очевидный факт
(иногда упоминаемый как свойство трактрисы прямой):\smallskip\\
{\em трактриса кривой,
построенная с поводком длины~$T$, является 
ортогональной траекторией семейства окружностей радиуса~$T$,
центры которых которых лежат на данной кривой}.

\Cfig{t}{.7\textwidth}{TrxOrtho}%
{Трактриса как ортогональная траектория семейства окружностей}

\Cfig{ht}{.7\textwidth}{TrxInv}{Ортогональные траектории семейств коцикличных окружностей}

%\Wfigr{12}{.3\textwidth}{TrxEnvel}{}
Рассмотрим трактрису окружности 
как ортогональную траекторию семейства 
%коцикличных 
окружностей радиуса~\GT{T},
центры  которых лежат на данной окружности радиуса~\GT{R}.
% (\Reffig{TrxEnvel}).
Это семейство имеет две огибающие~--- окружности радиусов~$R{+}T$ и ~$\abs{R{-}T}$.

На \Reffig[a]{TrxInv} сплошными линиями показано семейство
(полу)окружностей для случая $T\igt R$ и трактриса $AA_1$ (вида Т1)
в качестве ортогональной траектории семейства. 
Пунктиром показана инверсная конфигурация
\equa{%
    (x',y')=\frac{I^2}{x^2+y^2}\cdot(x,y),\qquad I^2=R^2-T^2
}
со степенью инверсии \LT{I^2}. Это {\em гиперболическая} инверсия: 
она, в отличие от обычной (эллиптической), предполагает 
дополнительно поворот на \Deg{180}. 
Образом точки $A$ является точка $B$, полуокружности переходят 
в свои дополнения до полной окружности.
Из сохранения углов при инверсии мы можем заключить, 
что одна ортогональная траектория преобразуется в другую,
и трактриса $AA_1$ (типа T1) преобразуется в реверсную трактрису $BB_1$.

На \Reffig[b]{TrxInv} \  $T\ilt R$.
Окружность инверсии, показанная штрих-пунктиром, 
в этом случае действительна: \GT{I^2}.
Получаем трактрисы T3 и T5 в качестве ортогональных траекторий,
также являющиеся инверсными образами друг друга.
Окружность инверсии совпадает с предельными окружностями
каждой из двух трактрис и 
является ещё одной~--- тривиальной~--- ортогональной траекторией~\eqref{TrxTriv}.

Наконец, при $T=R$ получаем семейство коцикличных окружностей,
показанное на \Reffig[c]{TrxInv}, и полярную трактрису Т2.
%в качестве ортогональной траектории семейства. 
При инверсии этой конфигурации относительно окружности~$A$
(огибающей семейства)
%(или любой другой окружности с центром в начале координат)
получим семейство прямых, касающихся~$A$, и их ортогональную
траекторию, являющуюся, как известно,
эвольвентой окружности~$A$ \cite[стр.\,252]{Savelov}.

Разбиение кругового кольца на \quo{криволинейные прямоугольники} с помощью
семейства окружностей и ортогонального ему семейства трактрис
показано на \Reffig[d]{TrxInv}.
%\smallskip

\newcommand{\Skip}[1]{}%
\Skip{%
Преобразуем уравнения трактрис T3,T5 и T1,T6 так, чтобы отразить
их инверсную связь. Трактрису T5 заставим вращаться в сторону увеличения полярного угла:
$\widetilde{\varphi}_5(t)\ieq\varphi_5(-t)$.
Заметим, что при соотношениях
\equa{%
   \begin{array}{ll}
      w_3=-w_5,\quad&
      \sqrt{\frac{1{-}w_3}{1{+}w_3}}\,t_3 = \sqrt{\frac{1{-}w_5}{1{+}w_5}}\,t_5,\\
      w_1=-w_6,   &
      \sqrt{\frac{w_1{-}1}{w_1{+}1}}\,t_1 = \sqrt{\frac{w_6{-}1}{w_6{+}1}}\,t_6\,.
   \end{array}
}
выполнены тождества 
\equa{%
       \varphi_3(t_3;w_3)\equiv\widetilde{\varphi}_5(t_5;w_5),\qquad
       \varphi_1(t_1;w_1)\equiv\varphi_6(t_6;w_6).
}
При замене параметра семейства~$w$ и параметра~$t$ самой кривой на
\equa{
    \omega=\sqrt{\abs{\frac{w{-}1}{w{+}1}}}\qquad
    \mbox{и}\qquad u\ieq\omega t
%    \qquad
%    \left(w_{1,6}\ieq\frac{1+\omega^2}{1-\omega^2},\quad
%          w_{3,5}\ieq\frac{1-\omega^2}{1+\omega^2}\right)
}
уравнения трактрис принимают вид 
\equa{%TrXY16}{%
 \begin{array}{ll}
   \begin{array}{l}
      0\ilt\omega_1\ilt1,\\
      1\ilt\omega_6\ilt\infty,\\
      0\ile u\ile\infty:
   \end{array}
   \quad   
   &
   \begin{array}{l}
       p_{1,6}(u)=\frac{2 R}{\abs{\omega^{-1}-\omega}}%
                  \sqrt{\frac{1+\omega^2u^2}{\omega^2+u^2}},\\
       \varphi_{1,6}(u)=(\omega^{-1}+\omega)\arctan u
               -\arctan(\omega u)-\arctan (\omega^{-1}u);
   \end{array}\\
    \,&\\      
   \begin{array}{l}
      0\ilt\omega_3\ilt1,\\
      1\ilt\omega_5\ilt\infty,\\
      0\ile u\ilt 1:
   \end{array}
   &   
   \begin{array}{l}
       p_{3,5}(u)=\frac{2 R}{\omega^{-1}+\omega}%
                  \sqrt{\frac{1+\omega^2u^2}{\omega^2+u^2}},\\
      \varphi_{3,5}(u)=\abs{\omega^{-1}-\omega}\atanh u%
               -\arctan\frac{\abs{\omega^{-1} -\omega} u}{1+u^2}.
   \end{array}
  \end{array}  
}
%
%Знаки модуля в $\varphi_{3,5}$ служат лишь
%для задания трактрисе T5 положительного направления вращения. 
Легко убедиться, что  при условиях $\omega_1\omega_6\ieq 1$,
$\omega_3\omega_5\ieq 1$
%(ранее $w_1\ieq{-w_6}$, $w_3\ieq{-w_5}$) 
выполнено 
\equa{%
    p_1(u;\omega_1){\cdot}p_6(u;\omega_6)\equiv T^2{-}R^2,\qquad
       p_3(u;\omega_3){\cdot}p_5(u;\omega_5)\equiv R^2{-}T^2.
}
}   %  Skip
\begin{center}
%\centerline
{$*\; *\; *$}
\end{center}

Отметим в заключение, что в последние годы 
в области Computer-Aided Design сильно возрос 
интерес к кривым с монотонно изменяющейся кривизной.
В рамках этих приложений их называют {\em спиралями} 
(об этом см. в \cite{PomiMain}).
Исследование спиралей и коллекционирование их образцов и привлекло
внимание автора к трактрисам окружности.

\section*{Литература}

%+Bibliography

%-Bibliography
\end{document}